\documentclass[]{interact}

\usepackage{epstopdf}
\usepackage[caption=false]{subfig}

\usepackage[numbers,sort&compress]{natbib}
\bibpunct[, ]{[}{]}{,}{n}{,}{,}
\renewcommand\bibfont{\fontsize{10}{12}\selectfont}
\makeatletter
\def\NAT@def@citea{\def\@citea{\NAT@separator}}
\makeatother

\theoremstyle{plain}
\newtheorem{theorem}{Theorem}[section]

\newtheorem{proposition}[theorem]{Proposition}

\theoremstyle{definition}

\theoremstyle{remark}

\newcommand{\ddu}{\ddot{\mbox{u}}}
\newcommand{\vect}{{\mbox{vec}}}

\usepackage{amssymb}
\usepackage{amssymb,epsfig,verbatim,amsmath,color}
\usepackage{hyperref}
\usepackage{lscape}
\usepackage{pdflscape}
\usepackage{afterpage}

\begin{document}


\title{Semiparametric Estimation for the Transformation Model with Length-Biased Data and Covariate Measurement Error}

\author{
\name{Li-Pang Chen\textsuperscript{a}\thanks{Correspondence author: Li-Pang Chen. 
Email: L358CHEN@uwaterloo.ca}}
\affil{\textsuperscript{a}Department of Statistics and Actuarial Science, 
University of Waterloo \\
200 University Ave W, Waterloo, ON N2L 3G1}
}

\maketitle

\begin{abstract}
Analysis of survival data with biased samples caused by left-truncation or length-biased sampling has received extensive interest. Many inference methods have been developed for various survival models. These methods, however, break down when survival data are typically error-contaminated. Although error-prone survival data commonly arise in practice, little work has been available in the literature for handling length-biased data with measurement error. In survival analysis, the transformation model is one of the frequently used models. However, methods of analyzing the transformation model with those complex features have not been fully explored. In this paper, we study this important problem and develop a valid inference method under the transformation model. We establish asymptotic results for the proposed estimators. The proposed method enjoys appealing features in that there is no need to specify the distribution of the covariates and the increasing function in the transformation model. Numerical studies are reported to assess the performance of the proposed method.
\end{abstract}

\begin{keywords}
Length-based sampling; measurement error; prevalent cohort; survival analysis; transformation model.
\end{keywords}

\section{Introduction} \label{Introduction}

In the study of a disease history, the time from the onset of an initiating event to the focused disease event (or failure) is usually the interest in the epidemiological and biomedical researches. One of the most attractive data comes from the prevalent sampling design, in which individuals only experience the initiating event but not the failure event before the recruiting time. Under this sampling scheme, individuals might not be observed because they experience the failure event before the recruiting time. Such a phenomenon caused by the delayed entry is called {\it left-truncation} and tends to produce a biased sample. Meanwhile, individuals who are recruited in the study may drop out or may not experience the failure event at the end of the study. It is called {\it right-censoring} in the dataset. To be specific, let $\xi$ be the calendar time of the recruitment and let $u$ and $r$ denote the calendar time of the initiating event (or the disease incidence) and the failure event, respectively, where $u<r$, and $u < \xi < r$.  Let $T^{*} = r-u$ be the failure time, and let $ A^{*} = \xi -u$ denote the truncation time. Assume that $A^\ast$ and $T^\ast$ are independent. Let $X^{*}$ be the associated covariate of dimension $p \times 1$.  Let $f(t)$ and $S(t)$ be the density function and the survivor function of the failure time $T^{*}$, respectively. Let $A,T$, and $X$ represent, respectively, the truncation time, the failure time, and the covariates for those subjects who are recruited in the study. That is, $(A,T,X)$ has the same joint distribution as $(A^\ast, T^\ast, X^\ast)$ given $T^\ast \geq A^\ast$. If $T^\ast < A^\ast$, then such an individual is not included in the study to contribute any information.
 
 Specifically, in the case of stable disease, the occurrence of disease onset follows the stationary Poisson process. To see this, we first impose two assumptions (e.g., Huang et al. 2012), including (a) the variables $(T^\ast, X^\ast)$ are independent of the disease incidence $u$, and (b) disease incidence occurs over calendar time at a constant rate. Based on these two assumptions, the joint density function of $\left(A^\ast, T^\ast \right)$ given $X^\ast = x$ is given by (Lancaster 1990, Chapter 3)
 \begin{eqnarray} \label{length-bias}
\frac{f(t|x)}{\mu(x)} I(t > a > 0), 
 \end{eqnarray}
 where $\mu(x) = \int_0^\infty \alpha f(\alpha | x) d \alpha$. Moreover, by (\ref{length-bias}), the failure time $T^\ast$ given $X^\ast = x$ has a length-biased density function $\frac{tf(t|x)}{\mu(x)}$.
  It implies that the truncation time $A^\ast$ follows the uniform distribution, and the survival time in the prevalent cohort has a length-biased sampling distribution since the probability of the survival time is proportional to the length of survival time. 
 
\begin{figure}[!h]
\centering{\includegraphics[scale=0.34]{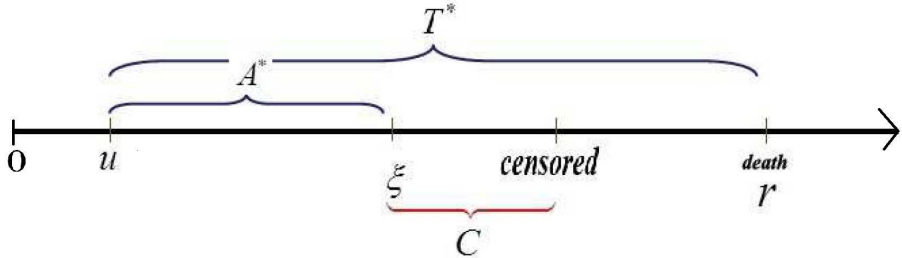}}
\caption{Schematic depiction of length-biased sampling.} \label{Calendar_time}
\end{figure}

  In addition, we let $C$ denote the censoring time for a recruited subject. Let $Y = \min \{ T, A+C \}$ be the observed survival time and let $\delta = I(T \leq A+C)$ be the indicator of a failure event. Figure~\ref{Calendar_time} gives an illustration of the relationship among those variables.

In the past literature, several models for the analysis of  biased samples were developed. To name a few, Huang and Qin (2013) and Chen (2018+) developed the estimation procedures for the additive hazards model for general left-truncated data where the distribution of the truncation time is unspecified. Qin and Shen (2010) and Huang et al. (2012) studied the Cox model for length-biased sampling. In addition, the transformation model is also a widely used model in survival analysis, which is formulated by
\begin{eqnarray} \label{AFT}
H\left( T^\ast\right) = -{X^\ast}^\top \beta + \epsilon,
\end{eqnarray}
where $H(\cdot)$ is an unknown increasing function, $\epsilon$ is a random variable with a known distribution, and $\beta$ is the $p \times 1$ vector of parameters of main interest. The transformation model gives a broad class of some frequently used models in survival analysis. Be more specific, when $\epsilon$ has an extreme value distribution, then $T^\ast$ follows the proportional hazards (PH) model; whereas when $\epsilon$ has a logistic distribution, then $T^\ast$ follows the proportional odds (PO) model. When biased samples occur, Shen et al. (2009) and Wang and Wang (2015) proposed valid methods to estimate the parameter of main interest based on length-biased sampling. In this paper, we mainly focus on the development of estimation method for the transformation model.

On the other hand, the other important feature in survival analysis is \textit{measurement error} in covariates. Assume that $X$ is an unobserved covariate and $W$ is a surrogate version of $X$. The classical additive model is frequently used, which is given by
\begin{eqnarray} \label{Meas-Model}
W = X + \eta,
\end{eqnarray}
where $\eta$ follows the normal distribution with mean zero and covariance matrix $\Sigma_\eta$. Assume that $X$ and $\eta$ are independent and $\Sigma_\eta$ is known for ease of discussion.

When the absence of left-truncation or length-biased sampling, there are many contributions on correcting the mismeasurement. For example, Nakamura (1992) proposed an approximate corrected partial likelihood method for the Cox model. Yi and Lawless (2007) developed the conditional expectation approach to correct the error effect for the Cox model. More discussions can be found in Yi (2017, Chapter 3).

However, little work has been available when those two features happen simultaneously. In the past literature, Chen and Yi (2018) proposed the corrected pseudo-likelihood estimation to estimate the parameter for the Cox model subject to left-truncation and right-censoring data and covariate measurement error. However, other survival models, such as the transformation model, have not been fully explored when those two complex features occur in the dataset. In this paper, our main purpose is to develop a valid estimation procedure to estimate the parameter for the transformation model with length-biased sampling and measurement error in covariate.

Our work is partially motivated by the Worcester Heart Attack Study (WHAS500) data (Hosmer et al. 2008) which involves biased and incomplete data. Basically, three types of time are recorded:  time of the hospital admission, time of the hospital discharge, and time of the last follow-up (which is either death or censoring time). The total follow-up length is defined as the time gap between the hospital admission and the last follow-up, and the hospital stay time is defined as the time length between the hospital admission and the hospital discharge. Data can only be collected for those individuals whose total follow-up length is larger than the hospital stay time, creating biased sample (e.g., Kalbfleisch and Prentice 2002, Section 1.3; Lawless 2003, Section 2.4).  It is interesting to study how the risk factors are associated with the survival times after the patients are discharged from the hospital. In addition, `blood pressure' and `body mass index' are variables in the data set, and they may be measured with error. To conduct sensible analyses, it is imperative to account for possible measurement error effects that are induced from error-prone covariates.  

The remainder is organized as follows. In Section~\ref{Main-Result}, we first present the proposed method to correct the error effect and derive the estimator. After that, we develop the theoretical result for the proposed method. Numerical results, including simulation study and real data analysis, are provided in Sections~\ref{Numerical} and~\ref{RDA}, respectively. Finally, we conclude the paper with discussions in Section~\ref{Summary}.


\section{Main Results} \label{Main-Result}
In this section, we propose a method to deal with the measurement error in covariates, and then derive the estimators of $\beta$ and $H(\cdot)$. The theoretical results of the estimators are also established in this section.

\subsection{Construction of The Estimating Function} \label{Method-1}
In this section, we review a method proposed by Wang and Wang (2015) and assume that the covariate $X_i$ is collected without mismeasurement. 

 Suppose we have a sample of $n$ subjects where for $i=1, \cdots, n$, $(Y_{i}, A_{i}, \delta_{i}, X_{i}, W_i)$ has the same distribution as $\left( Y,A,\delta,X ,W\right)$. Let $R_i(t) = I\left( A_i \leq t \leq Y_i \right)$ and $N_i(t) = I \left( Y_i \leq t, \delta_i = 1 \right)$ for $i = 1,\cdots,n$. Define 
\begin{eqnarray} \label{martingale}
M_i(t) = N_i(t) - \int_0^t R_i(u) \frac{\delta_i w(u)}{w(Y_i)} \lambda(u|X_i)
\end{eqnarray}
for $i = 1,\cdots,n$, where $\lambda(y|X_i)$ is the conditional hazard function of $Y_i$ given $X_i$, $w(t) = \int_0^t S_C(u)du$, and $S_C(t)$ is the survivor function of the variable $C$. It can be shown that (\ref{martingale}) is a square-integrable martingale (e.g., Anderson et al. 1993). Define $r(t,Y_i,\delta_i) = \frac{\delta_i w(t)}{w(Y_i)}$, and let $\widehat{r}(t,Y_i,\delta_i) = \frac{\delta_i \widehat{w}(t)}{\widehat{w}(Y_i)}$ denote the estimator of $r(t,Y_i,\delta_i)$, where $\widehat{w}(t) = \int_0^t \widehat{S}_C(u)du$ and $\widehat{S}_C(\cdot)$ is the Kaplan-Meier estimator. Therefore, based on (\ref{martingale}) and the counting process theory (e.g., Anderson et al. 1993), we can construct the following estimating equations
\begin{eqnarray*}
&\ & \sum \limits_{i=1}^n dM_i(t) = 0, \\
&\ & \sum \limits_{i=1}^n X_i dM_i(t) = 0,
\end{eqnarray*}
or equivalently,
\begin{eqnarray} \label{EE-X-1} 
&\ & \sum \limits_{i=1}^n \left[ dN_i(t) -  R_i(t) \widehat{r}(t,Y_i,\delta_i) d\Lambda\left\{X_i^\top \beta + H(t)\right\} \right] = 0,
\end{eqnarray}
and
\begin{eqnarray} \label{EE-X-2}
&\ & \sum \limits_{i=1}^n X_i \left[ dN_i(t) -  R_i(t) \widehat{r}(t,Y_i,\delta_i) d\Lambda\left\{X_i^\top \beta + H(t)\right\} \right] = 0,
\end{eqnarray}
where $\Lambda(\cdot)$ is the cumulative hazard function of $\epsilon$. When $\beta$ is fixed, let $\widehat{H}_X(t;\beta)$ denote the unique solution of (\ref{EE-X-1}). Next, replacing $H(\cdot)$ in (\ref{EE-X-2}) by $\widehat{H}_X(t;\beta)$ gives
\begin{eqnarray} \label{EE-X}
U_X(\beta) \triangleq  \sum \limits_{i=1}^n X_i \left[ dN_i(t) -  R_i(t) \widehat{r}(t,Y_i,\delta_i) d\Lambda\left\{X_i^\top \beta + \widehat{H}_X(t;\beta)\right\} \right] = 0.
\end{eqnarray}
Consequently, let $\widehat{\beta}_{true}$ denote the solution of $U_X(\beta)=0$ and it is the estimator of $\beta$ based on the covariate $X_i$. Moreover, $\widehat{H}_X(t;\widehat{\beta}_{true})$ is the estimator of $H(t)$.

Based on the true covariates $X_i$, we have the following proposition which is given in Theorem~1 of Wang and Wang (2015).
\begin{proposition}
Let $\beta_0$ be the true parameter. Under regularity conditions in  Wang and Wang (2015), $\widehat{\beta}_{true}$ is a consistent estimator and $\sqrt{n}\left( \widehat{\beta}_{true} - \beta_0 \right)$ converges weakly to a normal distribution with mean zero and variance-covariance matrix $\Sigma_\ast^{-1}\Sigma^\ast (\Sigma_\ast^{-1})^\top$, where $\Sigma_\ast$ and $\Sigma^\ast$ can be found in  Wang and Wang (2015).
\end{proposition}

However, as described in Section~\ref{Introduction}, we only observe $W_i$ instead of $X_i$. Therefore, replacing $X_i$ by $W_i$ in (\ref{EE-X}) gives the \textit{naive estimating function}
\begin{equation} \label{EE-W}
U_{naive}(\beta) \triangleq  \sum \limits_{i=1}^n W_i \left[ dN_i(t) -  R_i(t) \widehat{r}(t,Y_i,\delta_i) d\Lambda\left\{W_i^\top \beta + \widehat{H}_W(t;\beta)\right\} \right],
\end{equation}
where $\widehat{H}_W(t;\beta)$ is the solution of (\ref{EE-X-1}) replacing $X_i$ by $W_i$.
We let $\widehat{\beta}_{naive}$ denote the \textit{naive estimator} which is a solution of $U_{naive}(\beta) = 0$.

\subsection{Corrected Estimating Function} \label{Method-2}
It is well-known that the naive estimator $\widehat{\beta}_{naive}$ incurs the tremendous bias (e.g., Carroll et al. 2006; Yi 2017). To correct the mismeasurement and reduce the bias of the estimator, we employ the simulation-extrapolation (SIMEX) method (Cook and Stefaski 1994). The proposed procedure is the following three stages: 
\begin{description}
\item[Stage 1]Simulation\\
Let $B$ be a given positive integer and let $\mathcal{Z} = \left\{ \zeta_0,\zeta_1,\cdots,\zeta_M \right\}$ be a sequence of pre-specified values with $0 = \zeta_0<\zeta_1<\cdots<\zeta_M$. where $M$ is a positive integer, and $\zeta_M$ is pre-specified positive number such as $\zeta_M = 1$.

For a given subject $i$ with $i=1,\cdots,n$ and $b = 1,\cdots,B$, we generate $\eta_{b,i}$ from $N\left(0,\Sigma_\eta\right)$. Then for observed vector of covariates $W_i$, we define $W_i(b,\zeta)$ as
\begin{eqnarray} \label{SIMEX-1}
W_i(b,\zeta) = W_i + \sqrt{\zeta} \eta_{b,i}
\end{eqnarray}
for every $\zeta \in \mathcal{Z}$. Therefore, the conditional distribution of $W_i(b,\zeta)$ given $X_i$ is $N\left( X_i, (1+\zeta) \Sigma_\eta \right)$.

\item[Stage 2] Estimation\\
In this stage, replacing $X_i$ in (\ref{EE-X-1}) and (\ref{EE-X-2}) by $W_i(b,\zeta)$ defined in  (\ref{SIMEX-1}) gives new estimating equations
\begin{eqnarray} \label{EE-SIMEX-1} 
&\ & \sum \limits_{i=1}^n \left[ dN_i(t) -  R_i(t) \widehat{r}(t,Y_i,\delta_i) d\Lambda\left\{W_i^\top(b,\zeta) \beta + H(t)\right\} \right] = 0,
\end{eqnarray}
and
\begin{eqnarray} \label{EE-SIMEX-2}
&\ & \sum \limits_{i=1}^n W_i(b,\zeta) \left[ dN_i(t) -  R_i(t) \widehat{r}(t,Y_i,\delta_i) d\Lambda\left\{W_i^\top(b,\zeta) \beta + H(t)\right\} \right] = 0. \nonumber \\
\end{eqnarray}
Similar to the procedure in Section~\ref{Method-1}, we first fix $\beta$ and let $\widehat{H}(t;b,\zeta,\beta)$ denote the solution of (\ref{EE-SIMEX-1}). After that, replacing $H(\cdot)$ in (\ref{EE-SIMEX-2}) by $\widehat{H}(t;b,\zeta,\beta)$ gives the estimating equation 
\begin{eqnarray} \label{EE-SIMEX}
&\ \ &U_{SIMEX}(\beta) \nonumber \\
&\triangleq& \frac{1}{n} \sum \limits_{i=1}^n W_i(b,\zeta) \left[ dN_i(t) -  R_i(t) \widehat{r}(t,Y_i,\delta_i) d\Lambda\left\{W_i^\top(b,\zeta) \beta + \widehat{H}(t;b,\zeta,\beta)\right\} \right] \nonumber  \\
&=& 0.
\end{eqnarray}
for every $b = 1,\cdots,B$ and $\zeta \in \mathcal{Z}$. Let $\widehat{\beta}(b,\zeta)$ denote the solution of $U_{SIMEX}(\beta) = 0$. Moreover, we define
\begin{eqnarray} \label{SIMEX-2}
\widehat{\beta}(\zeta) = \frac{1}{B} \sum \limits_{b=1}^B \widehat{\beta}(b,\zeta).
\end{eqnarray}

\item[Stage 3] Extrapolation\\
By (\ref{SIMEX-2}), we have a sequence $\left\{ \left( \zeta, \widehat{\beta}(\zeta)\right) : \zeta \in \mathcal{Z} \right\}$. Then we fit a regression model to the sequence
\begin{eqnarray}
\widehat{\beta}(\mathcal{Z}) = \varphi\left(\mathcal{Z}, \Gamma\right) + \varrho,
\end{eqnarray}
where $\varphi(\cdot)$ is the user-specific regression function, $\Gamma$ is the associated parameter, and $\varrho$ is the noise term. The parameter $\Gamma$ can be estimated by the least square method, and we let $\widehat{\Gamma}$ denote the resulting estimate of $\Gamma$.

Finally, we calculate the predicted value
\begin{eqnarray} \label{SIMEX-3}
\widehat{\beta}_{SIMEX} = \varphi \left(-1, \widehat{\Gamma} \right)
\end{eqnarray}
and take $\widehat{\beta}_{SIMEX}$ as the \textit{SIMEX} estimator of $\beta$.

\item[Stage 4] Estimation of $H(\cdot)$ \\
Furthermore, we can also derive the estimator of the unknown function $H(\cdot)$. To do this, we first replace $\beta$ by $\widehat{\beta}_{SIMEX}$ in $\widehat{H}(t;b,\zeta,\beta)$, which gives $\widehat{H}(t;b,\zeta,\widehat{\beta}_{SIMEX})$. For every $t$ and $\zeta$, taking average with respect to $b$ gives $\widehat{H}(t;\zeta,\widehat{\beta}_{SIMEX}) = \frac{1}{B} \sum \limits_{b=1}^B \widehat{H}(t;b,\zeta,\widehat{\beta}_{SIMEX})$. Finally, similar to Stage 3 above, fitting a regression model and taking $\zeta = -1$ as a predicted value yields a final estimator $\widehat{H}(t;\widehat{\beta}_{SIMEX})$, also denoted as $\widehat{H}_{SIMEX}(t)$.
\end{description}

\subsection{Theoretical Results} \label{Theory}
In this section, we present the consistency and the asymptotic distribution of the proposed estimators $\widehat{\beta}_{SIMEX}$ and $\widehat{H}_{SIMEX}(\cdot)$, and the proofs are available in Appendix~\ref{pf}.

\begin{theorem} \label{Thm-SIMEX}
Let $\beta_0$ be the true parameter and let $H_0(t)$ denote the true increasing function. Under regularity conditions in Appendix~\ref{R.C.}, estimators $\widehat{\beta}_{SIMEX}$ and $\widehat{H}_{SIMEX}(\cdot)$ have the following properties:
\begin{itemize}
\item[(1)] $\widehat{\beta}_{SIMEX} \stackrel{p}{\longrightarrow} \beta_0$  as  $n \rightarrow \infty$;
\item[(2)] $\widehat{H}_{SIMEX}(t) \stackrel{p}{\longrightarrow} H_0(t)$  as  $n \rightarrow \infty$;
\item[(3)] $\sqrt{n} \left(\widehat{\beta}_{SIMEX} - \beta_0 \right) \stackrel{d}{\longrightarrow} N\left(0, \left\{\frac{\partial \varphi}{\partial \Gamma}\left(-1, \widehat{\Gamma} \right) \right\} \mathcal{Q} \left\{ \frac{\partial \varphi}{\partial \Gamma} \left(-1, \widehat{\Gamma} \right) \right\}^\top  \right)$  as  $n \rightarrow \infty$;
\item[(4)] $\sqrt{n} \left\{ \widehat{H}_{SIMEX}(t) - H_0(t) \right\}$ converges to the Gaussian process with mean zero and covariance function $E\left\{ \mathcal{H}_i(t) \mathcal{H}_i(s) \right\}$,
\end{itemize}
where the exact formulations of $\mathcal{Q}$ and $\mathcal{H}_i(t)$ are placed in Appendix~\ref{pf}.
\end{theorem}
\section{Numerical Study} \label{Numerical}

\subsection{Simulation Setup}
 We examine the setting where $\epsilon$ is generated from the extreme value distribution and the logistic distribution, and the truncation time $A^\ast$ is generated from the uniform distribution $U(0,1)$. Let $\beta = (\beta_{1}, \beta_{2})^\top$ denote a two-dimensional vector of parameters, and let $\beta_0 = (\beta_{10}, \beta_{20})^\top$ be the vector of true parameters where we set $\beta_0 = \left( 1,1 \right)^\top$. We consider a scenario where $X^\ast = (X_1^\ast, X_2^\ast)^\top$ are generated from a bivariate normal distribution with mean zero and variance-covariance matrix $\Sigma$, which is set as $ \left( \begin{array}{ c c} 
4 & 0.7 \\
0.7 & 3
\end{array}  
   \right)$. 
 Given $\epsilon$, $X^\ast$ and $\beta_0$, the failure time $T^\ast$ is generated from the model: 
\begin{eqnarray*} 
\log T^\ast = - \left( X_1^\ast \beta_{10} + X_2^\ast \beta_{20} \right) + \epsilon.
\end{eqnarray*}
Based on our two settings of $\epsilon$, the failure time $T^\ast$ follows the PH model and the PO model, respectively.
  Therefore, the observed data $(A,T,X)$ is collected from $(A^\ast,T^\ast,X^\ast)$ by conditioning on that $T^\ast \geq A^\ast$. We repeatedly generate data these steps we obtain a sample of a required size $n=200$. For the measurement error process, we consider model $(\ref{Meas-Model})$  with error $\eta \sim N \left( 0, \Sigma_\eta \right)$, where  $\Sigma_\eta$ is a diagonal matrix of $\sigma_\eta$ which is taken as $ 0.01$, $0.5$, and $0.75$, respectively.
  
We consider two censoring rates, say 25\% and 50\%, and let the  censoring time $C$ be generated from the uniform distribution $U(0,c)$, where $c$ is determined by a given censoring rate.
Consequently, $Y$ and $\delta$ are determined by $Y = \min \left\{ T, A+C \right\}$ and $\delta = I \left( T \leq A + C \right)$. In implementing the proposed method, we set $B = 500$ and partition the interval $[0,2]$ into subintervals with width $0.25$, and let the resulting cutpoints be the values of $\zeta$. We take the regression function $\varphi(\cdot)$ to be the quadratic polynomial function, which is a widely used function in many cases (e.g., Cook and Stefaski 1994; Carroll et al. 2006). Finally, 1000 simulations are run for each parameter setting. 

\subsection{Simulation Results}

We mainly examine three estimators which are discussed in Sections~\ref{Method-1} and~\ref{Method-2}, and they are labeled by Wang-Wang ($\widehat{\beta}_{true}$), Naive ($\widehat{\beta}_{naive}$), and Chen ($\widehat{\beta}_{SIMEX}$). We report the biases of estimates (Bias), the empirical variances (Var), the mean squared errors (MSE), and the coverage probabilities (CP) of those three estimators under the measurement error model (\ref{Meas-Model}). The results are reported in Table~\ref{tab:Sim}. 

First, the censoring rate and measurement degree have noticeable impact on each estimation methods.
As expected, biases and variance estimates increase as the censoring rate increases. 
When the measurement degree increases, biases of both $\widehat{\beta}_{naive}$ and $\widehat{\beta}_{SIMEX}$ are increasing, and the impact of the measurement error degrees seems more obvious on the naive estimator $\widehat{\beta}_{naive}$.

Within a setting with a given censoring rate and a measurement error degree, the naive method and the proposed method perform differently. When measurement error occurs, the performance of the proposed method is better than the naive method. The naive method produces considerable finite sample biases with coverage rates of 95\% confidence intervals significantly departing from the nominal level. The proposed method outputs satisfactory estimate with small finite sample biases and reasonable coverage rates of 95\% confidence intervals. Compared to the variance estimates produced by the naive approach, the proposed method which accounts for measurement error effects yield larger variance estimates, and this is the price paid to remove biases in point estimators. This phenomenon is typical in the literature of measurement error models. However, mean squared errors produced by the proposed method tends to be a lot smaller than those obtained from the naive method. Finally, we also present the numerical results of the estimator $\widehat{\beta}_{true}$ with the true covariate $X$. In general, $\widehat{\beta}_{true}$ gives the smallest bias and it is an efficient estimator with small variance. This result also verifies the validity of method proposed by Wang and Wang (2015).

\begin{landscape}
 \begin{table}
       \huge
     \caption{Numerical results for simulation study. `$\times$' denotes usage of the true covariate $X$; cr is the censoring rate; Bias is the difference between empirical mean and true value; Var is the empirical variance; MSE is the mean square error; and CP is the coverage probability.} \label{tab:Sim}

 \small

\center
   \renewcommand{\arraystretch}{0.25}
 \begin{tabular}{c c c c  ccccccccccccc}

 \\
 \hline\hline
model & cr  & $\sigma_{\eta}$  &  Method  & \multicolumn{4}{c} {Estimator of $\beta_1$ } & & \multicolumn{4}{c}{Estimator of $\beta_2$ } \\ \cline{5-8}  \cline{10-13}

 &  &  &  & Bias & Var & MSE & CP(\%) &  & Bias & Var  & MSE & CP(\%)
\\
 \hline 
PH & 25\% &  0.01 & Naive$(\widehat{\beta}_{naive})$ & -0.230 & 0.007 & 0.059 & 21.3 &  & -0.749 & 0.014 &  0.626 & 14.9\\
  &    &     & Chen$(\widehat{\beta}_{SIMEX})$    & 0.017 & 0.013  & 0.014 & 94.7 &  & 0.009 & 0.028 &  0.028 & 94.2   \\
      \\
  &    & 0.50 & Naive$(\widehat{\beta}_{naive})$ & -0.343 & 0.006  & 0.123 & 1.6 &  & -0.606 & 0.015 & 0.432 & 30.0\\
  &    &       & Chen$(\widehat{\beta}_{SIMEX})$        & 0.025 & 0.023 & 0.026 & 94.5 &  & 0.011 & 0.027  & 0.028 & 94.5\\
      \\
  &    & 0.75  & Naive$(\widehat{\beta}_{naive})$ &-0.347 & 0.005  & 0.125 & 0.3 &  & -0.636 & 0.016  & 0.465 & 23.8\\
  &    &       & Chen$(\widehat{\beta}_{SIMEX})$       & 0.025 & 0.023  & 0.023 & 94.8 &  & 0.019 & 0.025  & 0.025 & 93.9\\       
      \\
  &    &  $\times$ &  Wang-Wang$(\widehat{\beta}_{true})$  & 0.009 & 0.011 & 0.011 & 95.6 &  & 0.008 & 0.013  & 0.013 & 94.1\\ \cline{2-13}
& 50\% & 0.01 & Naive$(\widehat{\beta}_{naive})$ & -0.248 & 0.016  & 0.267 & 9.1 &  & -0.742 & 0.016  & 0.565 & 0.1\\
 &     &     & Chen$(\widehat{\beta}_{SIMEX})$       & 0.017 & 0.014  & 0.014 & 94.4 &  & 0.016 & 0.021  & 0.021 & 94.3\\
      \\
  &    & 0.50 & Naive$(\widehat{\beta}_{naive})$ &  -0.375 & 0.015  & 0.145 & 0.2 &  & -0.600 & 0.016  & 0.376 & 0.4\\
   &   &       & Chen$(\widehat{\beta}_{SIMEX})$        & 0.024 & 0.036  & 0.039 & 95.2 &  & 0.019 & 0.025 & 0.025 & 95.0\\
      \\
 &      & 0.75  & Naive$(\widehat{\beta}_{naive})$ & -0.360 & 0.014  & 0.134 & 0.1 &  & -0.630 & 0.014  & 0.413 & 0.2\\
 &     &       & Chen$(\widehat{\beta}_{SIMEX})$        &0.025 & 0.033 & 0.033 & 94.6 &  & 0.026 & 0.025  & 0.025 & 94.8\\   
       \\
  &    &  $\times$ &  Wang-Wang$(\widehat{\beta}_{true})$ & 0.008 & 0.013 & 0.013  & 95.3 &  & 0.006 & 0.014 & 0.014 & 95.7\\
  \hline 
PO & 25\% &  0.01 & Naive$(\widehat{\beta}_{naive})$ &  -0.250 & 0.009 &  0.072 & 23.0 & &  -0.729 & 0.015  & 0.557 & 0.4 \\
  &    &     & Chen$(\widehat{\beta}_{SIMEX})$    & 0.010 & 0.019  & 0.020 & 94.2 &  & 0.009 & 0.024  & 0.024 & 94.5   \\
      \\
  &    & 0.50 & Naive$(\widehat{\beta}_{naive})$ & -0.377 & 0.008  & 0.150 & 1.5 &  & -0.588 & 0.017  & 0.369 & 3.6\\
  &    &       & Chen$(\widehat{\beta}_{SIMEX})$        & 0.012 & 0.018  & 0.040 & 94.5 &  & 0.011 & 0.024 & 0.025 & 93.7\\
      \\
  &    & 0.75  & Naive$(\widehat{\beta}_{naive})$ &-0.362 & 0.007  & 0.138 & 1.1 &  & -0.619 & 0.014  & 0.405 & 1.4\\
  &    &       & Chen$(\widehat{\beta}_{SIMEX})$       & 0.016 & 0.018  & 0.018 & 94.3 &  & 0.015 & 0.022  & 0.022 & 94.6\\       
      \\
  &    &  $\times$ &  Wang-Wang$(\widehat{\beta}_{true})$  &0.005 & 0.015  & 0.015 & 94.6 &  & 0.005 & 0.013  & 0.014 & 94.7\\ \cline{2-13}
& 50\% & 0.01 & Naive$(\widehat{\beta}_{naive})$ & -0.268 & 0.016  & 0.273 & 10.1 &  & -0.842 & 0.016  & 0.574 & 1.3\\
 &     &     & Chen$(\widehat{\beta}_{SIMEX})$       & 0.016 & 0.024  & 0.024 & 94.6 &  & 0.016 & 0.027  & 0.027 & 94.5\\
      \\
  &    & 0.50 & Naive$(\widehat{\beta}_{naive})$ &  -0.388 & 0.016  & 0.168 & 1.4 &  & -0.600 & 0.016  & 0.376 & 1.4\\
   &   &       & Chen$(\widehat{\beta}_{SIMEX})$        & 0.027 & 0.036  & 0.037 & 94.2 &  & 0.021 & 0.026 & 0.026 & 95.1\\
      \\
 &      & 0.75  & Naive$(\widehat{\beta}_{naive})$ & -0.410 & 0.017  & 0.185 & 1.9 &  & -0.630 & 0.018  & 0.413 & 1.2\\
 &     &       & Chen$(\widehat{\beta}_{SIMEX})$        &0.028 & 0.036 & 0.036 & 94.6 &  & 0.025 & 0.027  & 0.027 & 94.6\\   
       \\
  &    &  $\times$ &  Wang-Wang$(\widehat{\beta}_{true})$ & 0.010 & 0.013 & 0.013  & 95.1 &  & 0.008 & 0.014 & 0.014 & 95.4\\
 
 \hline\hline
 \\
\end{tabular}
\end{table}
\end{landscape}

 
\section{Data Analysis} \label{RDA}
 
In this section, we apply the proposed method to analyze the data arising from the Worcester Heart Attack Study (WHAS500), which are described in Section~\ref{Introduction}. 
Discussed by Hosmer et al. (2008), a survival time was defined as the time since a subject  was admitted to the hospital. We are interested in studying survival times of patients who were discharged alive from the hospital. Hence, a selection criterion was imposed that only those subjects who were discharged alive were eligible to be included in the analysis. That is, individuals were not enrolled in the analysis if they died before discharging from the hospital, hence biased sample occurs. With such a criterion, a sample of size 461 was available.
In this data set, the censoring rate is 61.8\%. 
To be more specific, the total length of follow-up (lenfol) is the survival time (i.e., $Y_i = \min \left(T_i,C_i \right)$), the length of hospital stay (los) is the truncation time (i.e., $A_i$), and the vital status at last follow-up (fstat) is $\delta_i$. In our analysis, the covariates include the body mass index (BMI) and the blood pressure (BP) of a patient. Suppose that BMI and BP are subject to measurement error, we let $W = \left( W_1, W_2 \right)^\top$ where $W_1$ and $W_2$ denote BMI and BP, respectively, and consider the measurement error model $(\ref{Meas-Model})$. Let $\beta = \left( \beta_1,\beta_2 \right)^\top$, where $\beta_k$ is the associated parameter with $W_k$ for $k=1,2$. We fit the regression model (\ref{AFT}) with $\epsilon$ being the extreme value distribution and the logistic distribution, respectively.

In this data set, there is no additional data source, such as a validation subsample or replicated measurements which is often required to describe the measurement error process (e.g., Carroll et al. 2006; Yi 2017). To get around this and understand the impact of measurement error on estimation, we carry out sensitivity analyses.  Specifically, let $\Sigma$ be the sample covariance matrix of $W$, and for sensitivity analyses we consider $\Sigma_\eta = \Sigma + \Sigma_e$ to be the covariance matrix for the measurement error model (\ref{Meas-Model}), where $\Sigma_e$ is the diagonal matrix with diagonal elements being a common value $\sigma_e \in [0,1]$.

In Table~\ref{Table: RDA}, we report the point estimates (EST), the standard errors and p-values for the estimators $\widehat{\beta}_{SIMEX}$ for the cases with $\sigma_e = 0.15$, $0.50$ and $0.75$, respectively, corresponding to minor, moderate and large measurement error, and report the results of $\widehat{\beta}_{naive}$ by directly using $W$ in the analysis. All the point estimates produced by our proposed method are fairly close as observed, and the results are fairly stable regardless of the degree of measurement error. Both our proposed method and the naive estimator suggest that two covariates are significant regardless of the specification of the distribution of $\epsilon$. 

\begin{table}
 \caption{Sensitivity analyses result of Worcester Heart Attack Study Data} \label{Table: RDA}
 \centering
 
 \begin{tabular}{cc c  c ccccccc}
 \\
 \hline\hline
&  &   &  \multicolumn{3}{c}{ Estimator $\beta_1$} & & \multicolumn{3}{c}{Estimator $\beta_2$}
\\ \cline{4-6} \cline{8-10} 
 Model & $\sigma_e$ & method  & EST & S.E & p-value & & EST & S.E & p-value 
\\ [0.5 ex]
PH & 0.15   &$\widehat{\beta}_{SIMEX}$ &-3.297 & 0.011 & 0.000 & & -1.448 & 0.174 & 8.661e-17 \\
  

&0.50   &$\widehat{\beta}_{SIMEX}$ &-1.370 & 0.011 & 0.000 & & -1.405 & 0.179 & 4.188e-15 \\
   

& 0.75   &$\widehat{\beta}_{SIMEX}$ & -2.091 & 0.016 & 0.000 & & -1.382 & 0.152 &  9.714e-20 \\

&$\times$   &$\widehat{\beta}_{naive}$ &-1.694 &  0.010 & 0.000 & & -0.834  & 0.109 & 1.988e-14 \\
\hline
PO & 0.15   &$\widehat{\beta}_{SIMEX}$ &-1.538 & 0.015 & 1.169e-22 & & -1.570 & 0.174 & 6.765e-100 \\
  

&0.50   &$\widehat{\beta}_{SIMEX}$ &-1.444 & 0.012 &  2.374e-33 & & -1.250 & 0.174 & 2.816e-33 \\
   

& 0.75   &$\widehat{\beta}_{SIMEX}$ & -1.933 & 0.012 & 1.540e-56 & & -1.102 & 0.162 & 1.120e-70 \\

&$\times$   &$\widehat{\beta}_{naive}$ &-0.894 & 0.011 & 4.391e-16 & & -0.326 & 0.107 & 4.029e-12 \\
\hline\hline
\end{tabular}
\end{table}

\section{Discussion} \label{Summary}

In this article, we focus the discussion on the transformation model with length-biased sampling and develop a valid method to correct the covariate measurement error and derive an efficient estimator. In this article, we also establish the large sample properties, and the numerical results guarantee that our proposed method outperforms. Although we only focus on the simple structure of the measurement error model, our method can easily be extended to complex measurement error models or additional information, such as repeated measurement or validation data. In addition, there are still many challenges in this topic, such as analysis without specifying the distribution of the truncation time $A^\ast$ or the discussion of time-dependent covariates with mismeasurement. These topics are also our researches in the future.

\appendix

\section{Regularity Conditions} \label{R.C.}
\begin{itemize}
\item[(C1)] $\Theta$ is a compact set, and the true parameter value $\beta_0$ is an interior point of $\Theta$.
\item[(C2)] Let $\tau$ be the finite maximum support of the failure time.
\item[(C3)] The $\left\{ A_i, Y_i,X_i \right\}$ are independent and identically distributed for $i=1,\cdots,n$. 
\item[(C4)] The covariate $X_i$ is bounded.
\item[(C5)] Conditional on $X_i^\ast$, $\left( T_i^\ast,  X_i^\ast \right)$ are independent of $A_i^\ast$.
\item[(C6)] Censoring time $C_i$ is non-informative. That is, the failure time $T_i$ and the censoring time $C_i$ are independent, given the covariates $X_i$.
\item[(C7)] The regression function $\varphi(\cdot)$ is true, and its first order derivative exists.
\end{itemize}
Condition (C1) is a basic condition that is used to derive the maximizer of the target function. (C2) to (C6) are standard conditions for survival analysis, which allow us to obtain the sum of i.i.d. random variables and hence to derive the asymptotic properties of the estimators. Condition (C7) is a common assumption in SIMEX method.

\setcounter{equation}{0}
\renewcommand\theequation{B.\arabic{equation}}
\section{Proof of Theorem~\ref{Thm-SIMEX}} \label{pf}
Before presenting the proof, we first define some notation. For any $s,t \in [0,\tau]$, $b = 1,\cdots,B$, and $\zeta \in \mathcal{Z}$, let
\begin{eqnarray*}
B_1(t;b,\zeta) &=& E \left[ R_i(t) r(t,Y_i,\delta_i) \lambda'\left\{ W_i^\top(b,\zeta) \beta(b,\zeta) + H_0(t) \right\} \right], \\
B_2(t;b,\zeta) &=& E \left[ R_i(t) r(t,Y_i,\delta_i) \lambda\left\{ W_i^\top(b,\zeta) \beta(b,\zeta)  + H_0(t) \right\} \right], \\
B_1^W(t;b,\zeta) &=& E \left[ W_i(b,\zeta) R_i(t) r(t,Y_i,\delta_i) \lambda'\left\{ W_i^\top(b,\zeta) \beta(b,\zeta)  + H_0(t) \right\} \right], \\
B_2^W(t;b,\zeta) &=& E \left[ W_i(b,\zeta) R_i(t) r(t,Y_i,\delta_i) \lambda\left\{ W_i^\top(b,\zeta) \beta(b,\zeta)  + H_0(t) \right\} \right], \\
B(t,s;b,\zeta) &=& \exp \left\{ \int_s^t B_2^{-1}(u;b,\zeta) B_1(u;b,\zeta) dH_0(u)  \right\},
\end{eqnarray*}
where $\lambda'(\cdot)$ is the derivative of $\lambda(\cdot)$. We further define
\begin{eqnarray*}
\mathbf{w}(t) &=& B_2^{-1}(t;b,\zeta) \left[ B_2^W(t;b,\zeta) \right. \\
& \ \ & \left. + \int_t^\tau \left\{ B_1^W(s;b,\zeta) - B_2^{-1}(s;b,\zeta) B_2^W(s;b,\zeta) B_1(s;b,\zeta) \right\} B(t,s;b,\zeta) dH_0(s) \right], \\
a(t) &=& \frac{1}{\pi(t)} \int_0^\tau E \left[ \frac{W_i(b,\zeta) \delta_i R_i(s)}{w(Y_i)} \left( I(t \leq s) \int_t^s S_C(u) du \right. \right. \\
&\ \ & \left. \left. - \frac{w(s)}{w(Y_i)} I(t \leq Y_i ) \int_t^{Y_i} S_C(u) du\right) d\Lambda\left\{ W_i^\top (b,\zeta) \beta(b,\zeta) + H_0(s) \right\} \right],
\end{eqnarray*}
where $\pi(t) = P\left( Y_i - A_i \geq t \right)$.
\\
\\
\underline{Proof of Theorem~\ref{Thm-SIMEX} (1)}:\\
The Kaplan-Meier estimator $\widehat{S}_C(t)$ over $[0,\tau]$ is uniformly consistency to $S_C(t)$ in the sense that $\sup \limits_{t \in [0,\tau]} \left| \widehat{S}_C(t) - S_C(t) \right| \stackrel{p}{\longrightarrow} 0$ as $n \rightarrow \infty$ (Pollard 1990; van der Vaart 1998). It implies that as $n \rightarrow \infty$,
\begin{eqnarray}\label{pf-1-1}
\sup \limits_{t \in [0,\tau]} \left| \widehat{w}(t) - w(t) \right| \stackrel{p}{\longrightarrow} 0
\end{eqnarray}
and
\begin{eqnarray}\label{pf-1-2}
\sup \limits_{t \in [0,\tau]} \left| \widehat{r}(t,Y_i,\delta_i) - r(t,Y_i,\delta_i) \right| \stackrel{p}{\longrightarrow} 0
\end{eqnarray}
since $\sup \limits_{t \in [0,\tau]} \left| \int_0^t \widehat{S}_C(u)du - \int_0^t S_C(u)du \right| \leq \sup \limits_{t \in [0,\tau]}  \int_0^t\left| \widehat{S}_C(u)- S_C(u)\right|du  \stackrel{p}{\longrightarrow} 0$ as $n \rightarrow \infty$. Furthermore, by the similar derivation in Step 1 of Wang and Wang (2015), we have
\begin{eqnarray} \label{pf-1-2-1}
\left.\frac{\partial \widehat{H}(t;b,\zeta,\beta)}{\partial \beta} \right|_{\beta = \beta_0} &=& - \int_0^t \frac{B(s,t;b,\zeta)}{B_2(s;b,\zeta)} B_1^W (s;b,\zeta) dH_0(s) + o_p(1) \nonumber \\
&\triangleq& A(t) + o_p(1).
\end{eqnarray}

Let $\beta(b,\zeta)$ be a solution of $E\left\{U_{SIMEX}(\beta)\right\} = 0$. Since $\widehat{\beta}(b,\zeta)$ is a solution of $U_{SIMEX}(\beta) = 0$. By (\ref{pf-1-1}), (\ref{pf-1-2}), (\ref{pf-1-2-1}), and the Uniformly Law of Large Numbers (van der Vaart 1998), we have that $\frac{1}{n} U_{SIMEX}(\beta)$ converges uniformly to $E\left\{U_{SIMEX}(\beta)\right\}$. Then we have that as $n \rightarrow \infty$,
\begin{eqnarray} \label{pf-1-3}
\widehat{\beta}(b,\zeta) \stackrel{p}{\longrightarrow} \beta(b,\zeta).
\end{eqnarray}
By definition (\ref{SIMEX-2}), taking averaging with respect to $b$ on both sides of (\ref{pf-1-3}) gives that as $n \rightarrow \infty$,
\begin{eqnarray}  \label{pf-1-4}
\widehat{\beta}(\zeta) \stackrel{p}{\longrightarrow} \beta(\zeta)
\end{eqnarray}
for every $\zeta \in \mathcal{Z}$.
By (\ref{pf-1-4}), we can show that as $n \rightarrow \infty$,
\begin{eqnarray} 
\widehat{\Gamma} \stackrel{p}{\longrightarrow} \Gamma.
\end{eqnarray}
Since $\widehat{\beta}_{SIMEX} = \varphi\left(-1,\widehat{\Gamma} \right)$, therefore, by the continuous mapping theorem, we have that as $n \rightarrow \infty$,
\begin{eqnarray} \label{pf-1-5}
\widehat{\beta}_{SIMEX} \stackrel{p}{\longrightarrow} \beta_0.
\end{eqnarray}
\ \\
\underline{Proof of Theorem~\ref{Thm-SIMEX} (2)}:\\
By (\ref{pf-1-5}), we have $\widehat{H}(t;b,\zeta,\widehat{\beta}_{SIMEX}) - \widehat{H}(t;b,\zeta,\beta_0) = o_p(1)$ for every $t \in [0,\tau]$, b, and $\zeta$. Taking average with respect to $b$ gives $\widehat{H}(t;\zeta,\widehat{\beta}_{SIMEX}) - \widehat{H}(t;\zeta,\beta_0) = o_p(1)$. On the other hand, by the Uniformly Law of Large Numbers and similar derivations in Wang and Wang (2015) with $\zeta \rightarrow -1$, we have that as $n \rightarrow \infty$, $\widehat{H}(t;\beta_0) - H_0(t) \stackrel{p}{\longrightarrow} 0$ for all $t \in [0,\tau]$. Therefore, we conclude that as $n \rightarrow \infty$ and $\zeta \rightarrow -1$, $\widehat{H}_{SIMEX}(t) - H_0(t) \stackrel{p}{\longrightarrow} 0$ by the fact that $\widehat{H}(t;-1,\widehat{\beta}_{SIMEX}) - H_0(t) = \widehat{H}(t;-1,\widehat{\beta}_{SIMEX}) - \widehat{H}(t;-1,\beta_0) + \widehat{H}(t;-1,\beta_0) - H_0(t)$.
\\
\\
\underline{Proof of Theorem~\ref{Thm-SIMEX} (3)}:\\
For $b = 1,\cdots,B$ and $\zeta \in \mathcal{Z}$, applying the Taylor series expansion on (\ref{EE-SIMEX}) around ${\beta}(b,\zeta)$ gives
\begin{eqnarray*}
0 &=& U_{SIMEX}\left( \widehat{\beta}(b,\zeta) \right) \\
&=& U_{SIMEX}\left( {\beta}(b,\zeta) \right) + \frac{\partial U_{SIMEX}\left( {\beta}(b,\zeta) \right)}{\partial \beta} \left\{ \widehat{\beta}(b,\zeta) - {\beta}(b,\zeta) \right\} + o_p\left( \frac{1}{\sqrt{n}}\right),
\end{eqnarray*}
or equivalently,
\begin{eqnarray} \label{pf-3-1}
\sqrt{n} \left\{ \widehat{\beta}(b,\zeta) - {\beta}(b,\zeta) \right\} &=& \left(-\frac{\partial U_{SIMEX}\left( {\beta}(b,\zeta) \right)}{\partial \beta} \right)^{-1} \sqrt{n} U_{SIMEX}\left( {\beta}(b,\zeta) \right) \nonumber \\
&\ \ & + o_p\left( 1 \right).
\end{eqnarray}

By (\ref{pf-1-1}), (\ref{pf-1-2}), (\ref{pf-1-2-1}), and the Uniformly Law of Large Numbers, we have that as $n \rightarrow \infty$, \begin{eqnarray} \label{pf-3-2}
\left(-\frac{\partial U_{SIMEX}\left( {\beta}(b,\zeta) \right)}{\partial \beta} \right) \stackrel{p}{\longrightarrow} \mathcal{A}\left(b,\zeta \right),
\end{eqnarray}
where 
\begin{eqnarray*}
\mathcal{A}\left(b,\zeta \right) &=& - E\left[ \int_0^\tau \left\{ W_i(b,\zeta) - \mathbf{w}(t) \right\} W_i^\top(b,\zeta) R_i(t) \right. \\
&\ \ & \times \left. r(t,Y_i,\delta_i) \lambda' \left\{ W_i^\top (b,\zeta) \beta (b,\zeta) + H_0(t) \right\} dH_0(t) \right].
\end{eqnarray*}

On the other hand, by the similar derivations in Wang and Wang (2015) and standard algebra, we can derive $U_{SIMEX}({\beta}(b,\zeta))$ as a sum of i.i.d. random functions. The exact formulation is given by
\begin{eqnarray}  \label{pf-3-3}
\sqrt{n} U_{SIMEX}({\beta}(b,\zeta)) = \frac{1}{\sqrt{n}} \sum \limits_{i=1}^n \Psi_i(b,\zeta) + o_p(1),
\end{eqnarray}
where
\begin{eqnarray*}
\Psi_i(b,\zeta) &=& \int_0^\tau \left[ \left\{W_i(b,\zeta) - \mathbf{w}(t)\right\} dM_i(t) a(t) dM_{C}(t) \right],
\end{eqnarray*}
$M_{C} = I\left(Y_i \leq t, \delta_i = 0 \right) - \int_0^t I\left(Y_i \geq u \right) d\Lambda_C(u)$, and $\Lambda_C(\cdot)$ is the cumulative hazards function of $C$.

Combining (\ref{pf-3-3}) and (\ref{pf-3-2}) with (\ref{pf-3-1}) yields
\begin{eqnarray} \label{pf-3-4}
\sqrt{n} \left\{ \widehat{\beta}(b,\zeta) - {\beta}(b,\zeta) \right\} = \frac{1}{\sqrt{n}} \sum \limits_{i=1}^n \mathcal{A}^{-1}\left(b,\zeta \right) \Psi_i(b,\zeta) + o_p\left( 1 \right).
\end{eqnarray}
By (\ref{SIMEX-2}), taking average with respect to $b$ on both sides of (\ref{pf-3-4}) gives
\begin{eqnarray} \label{pf-3-5}
\sqrt{n} \left\{ \widehat{\beta}(\zeta) - {\beta}(\zeta) \right\} = \frac{1}{\sqrt{n}} \sum \limits_{i=1}^n \Phi_i(\zeta) + o_p\left( 1 \right)
\end{eqnarray}
for $\zeta \in \mathcal{Z}$, where $ \Phi_i(\zeta) =  \frac{1}{B} \sum \limits_{b=1}^B \mathcal{A}^{-1}\left(b,\zeta \right) \Psi_i(b,\zeta)$.

Let $\widehat{\beta}(\mathcal{Z}) = \vect\left\{\widehat{\beta}(\zeta) : \zeta \in \mathcal{Z}\right\}$ denote the vectorization of estimator $\widehat{\beta}(\zeta)$ with every $\zeta \in \mathcal{Z}$. By the Central Limit Theorem on (\ref{pf-3-5}), we have that as $n \rightarrow \infty$,
\begin{eqnarray} \label{pf-3-6}
\sqrt{n} \left\{ \widehat{\beta}(\mathcal{Z}) - {\beta}(\mathcal{Z}) \right\} \stackrel{d}{\longrightarrow} N\left(0, \Omega\left( \mathcal{Z} \right) \right),
\end{eqnarray} 
where $\Omega\left( \mathcal{Z} \right) = \text{cov}\left\{ \Phi_i(\mathcal{Z}) \right\}$. By the Taylor series expansion on $\varphi\left( \mathcal{Z}, \Gamma \right)$ with respect to $\Gamma$, we have
\begin{eqnarray} \label{pf-3-7}
\varphi\left( \mathcal{Z}, \widehat{\Gamma} \right) - \varphi\left( \mathcal{Z}, \Gamma \right) \approx \frac{\partial \varphi\left( \mathcal{Z}, \Gamma \right)}{\partial \Gamma} \left(\widehat{\Gamma} - \Gamma \right).
\end{eqnarray}
Let $\mathcal{C} = \frac{\partial \varphi\left( \mathcal{Z}, \Gamma \right)}{\partial \Gamma}$ and $\mathcal{D} = \left\{\frac{\partial \varphi\left( \mathcal{Z}, \Gamma \right)}{\partial \Gamma}\right\}^\top \frac{\partial \varphi\left( \mathcal{Z}, \Gamma \right)}{\partial \Gamma}$. Combining (\ref{pf-3-6}) and (\ref{pf-3-7}) gives that as $n \rightarrow \infty$,
\begin{eqnarray} \label{pf-3-8}
\sqrt{n} \left(\widehat{\Gamma} - \Gamma \right) \stackrel{d}{\longrightarrow} N\left(0, \mathcal{D}^{-1} \mathcal{C} \Omega(\mathcal{Z}) \mathcal{C}^\top \mathcal{D}^{-1}  \right).
 \end{eqnarray}
 Finally, let $\mathcal{Q} = \mathcal{D}^{-1} \mathcal{C} \Omega(\mathcal{Z}) \mathcal{C}^\top \mathcal{D}^{-1}$. Since the SIMEX estimator is defined by $\widehat{\beta}_{SIMEX} = \varphi\left(-1, \widehat{\Gamma} \right)$, then combining (\ref{pf-3-7}) and (\ref{pf-3-8}) with $\zeta \rightarrow -1$ and applying the delta method give that as $n \rightarrow \infty$,
 \begin{eqnarray*} 
\sqrt{n} \left(\widehat{\beta}_{SIMEX} - \beta_0 \right) \stackrel{d}{\longrightarrow} N\left(0, \left\{\frac{\partial \varphi}{\partial \Gamma}\left(-1, \widehat{\Gamma} \right) \right\} \mathcal{Q} \left\{ \frac{\partial \varphi}{\partial \Gamma} \left(-1, \widehat{\Gamma} \right) \right\}^\top  \right).
 \end{eqnarray*}
 \ \\ 
 \underline{Proof of Theorem~\ref{Thm-SIMEX} (4)}:\\
 We first consider the expression of $\sqrt{n} \left\{ \widehat{H}(t;b,\zeta,\widehat{\beta}_{SIMEX}) - \widehat{H}(t;b,\zeta,\beta_0) \right\}$. By the Taylor series expansion with respect to $\beta$, we have
 \begin{eqnarray} \label{pf-4-1}
&& \sqrt{n} \left\{ \widehat{H}(t;b,\zeta,\widehat{\beta}_{SIMEX}) - \widehat{H}(t;b,\zeta,\beta_0) \right\} \nonumber \\
&=& \frac{\partial \widehat{H}(t;b,\zeta,\beta_0)}{\partial \beta} \sqrt{n} \left( \widehat{\beta}_{SIMEX} - \beta_0 \right) \nonumber \\ 
&=& A(t) \sqrt{n} \left( \widehat{\beta}_{SIMEX} - \beta_0 \right) + o_p(1) \nonumber \\
&=& A(t) \frac{1}{\sqrt{n}} \sum \limits_{i=1}^n \left\{\frac{\partial \varphi}{\partial \Gamma}\left(-1, \widehat{\Gamma} \right) \right\} \mathcal{D}^{-1} \mathcal{C} \Phi_i(\mathcal{Z}) + o_p(1),
 \end{eqnarray}
where the second equality is due to (\ref{pf-1-2-1}) and the third equality is due to (\ref{pf-3-8}).

By similar derivations in Step 3 of Wang and Wang (2015), we have
\begin{eqnarray} \label{pf-4-2}
\sqrt{n} \left\{  \widehat{H}(t;b,\zeta,\beta_0) - H_0(t) \right\} = \frac{1}{\sqrt{n}} \sum \limits_{i=1}^n \int_0^t \frac{B(s,t;b,\zeta)}{B_2(s;b,\zeta)} dM_i(s) + o_p(1).
\end{eqnarray}
Then combining (\ref{pf-4-1}) and (\ref{pf-4-2}) gives
\begin{eqnarray} \label{pf-4-3}
\sqrt{n} \left\{  \widehat{H}(t;b,\zeta,\widehat{\beta}_{SIMEX}) - H_0(t) \right\} = \frac{1}{\sqrt{n}} \sum \limits_{i=1}^n \mathcal{T}_i(t;b,\zeta) + o_p(1),
\end{eqnarray}
where $\mathcal{T}_i(t;b,\zeta) = A(t) \left\{\frac{\partial \varphi}{\partial \Gamma}\left(-1, \widehat{\Gamma} \right) \right\} \mathcal{D}^{-1} \mathcal{C} \Phi_i(\mathcal{Z}) + \int_0^t \frac{B(s,t;b,\zeta)}{B_2(s;b,\zeta)} dM_i(s)$. Taking average on both sides of (\ref{pf-4-3}) with respect to $b$ yields
\begin{eqnarray} \label{pf-4-4}
\sqrt{n} \left\{  \widehat{H}(t;\zeta,\widehat{\beta}_{SIMEX}) - H_0(t) \right\} = \frac{1}{\sqrt{n}} \sum \limits_{i=1}^n \mathcal{T}_i(t;\zeta) + o_p(1),
\end{eqnarray}
where $\mathcal{T}_i(t;\zeta) = \frac{1}{B} \sum \limits_{b=1}^B \mathcal{T}_i(t;b,\zeta)$.

Suppose that $\varphi_H(\zeta,\Gamma_H(t))$ is a function with the same conditions in (C7), and $\Gamma_H(t)$ is the associated parameter depending on time $t$. For $t \in [0,\tau]$ and $\zeta \in \mathcal{Z}$, we fit a regression model on $\widehat{H}(t;\zeta,\widehat{\beta}_{SIMEX})$ and $\varphi_H(\zeta,\Gamma_H(t))$, and derive the estimator of $\Gamma_H(t)$ which is denoted by $\widehat{\Gamma}_H(t)$. Furthermore, similar to the derivations in (\ref{pf-3-7}), we have  
\begin{eqnarray} \label{pf-4-5}
\varphi_H\left( \mathcal{Z}, \widehat{\Gamma}_H(t) \right) - \varphi_H\left( \mathcal{Z}, \Gamma_H(t) \right) \approx \frac{\partial \varphi_H\left( \mathcal{Z}, \Gamma_H(t) \right)}{\partial \Gamma_H(t)} \left\{ \widehat{\Gamma}_H(t) - \Gamma_H(t) \right\}. 
\end{eqnarray}
Let $\mathcal{U}(t) = \frac{\partial \varphi_H\left( \mathcal{Z}, \Gamma_H(t) \right)}{\partial \Gamma_H(t)}$ and $\mathcal{V}(t) = \left\{\frac{\partial \varphi_H\left( \mathcal{Z}, \Gamma_H(t) \right)}{\partial \Gamma_H(t)}\right\}^\top \frac{\partial \varphi_H\left( \mathcal{Z}, \Gamma_H(t) \right)}{\partial \Gamma_H(t)}$. Combining (\ref{pf-4-4}) and (\ref{pf-4-5}) yields
\begin{eqnarray} \label{pf-4-6}
\sqrt{n} \left\{ \widehat{\Gamma}_H(t) - \Gamma_H(t) \right\} = \frac{1}{\sqrt{n}} \sum \limits_{i=1}^n \mathcal{V}^{-1}(t) \mathcal{U}(t) \mathcal{T}_i(t;\mathcal{Z}) + o_p(1),
\end{eqnarray}
and since the estimator $\widehat{H}_{SIMEX}(t)$ is a predicted value of $\varphi_H(\zeta,\widehat{\Gamma}_H(t))$ by taking $\zeta \rightarrow -1$, then by (\ref{pf-4-6}) and the delta method with $\zeta \rightarrow -1$, we obtain
\begin{eqnarray} \label{pf-4-7}
& & \sqrt{n} \left\{ \widehat{H}_{SIMEX}(t) - H_0(t) \right\} \nonumber \\
&=& \frac{1}{\sqrt{n}} \sum \limits_{i=1}^n \left\{ \frac{\partial \varphi_H(-1,\widehat{\Gamma}_H(t))}{\partial \Gamma_H(t)} \right\} \mathcal{V}^{-1}(t) \mathcal{U}(t) \mathcal{T}_i(t;\mathcal{Z}) + o_p(1) \nonumber \\
&\triangleq& \frac{1}{\sqrt{n}} \sum \limits_{i=1}^n \mathcal{H}_i(t) + o_p(1).
\end{eqnarray}
where $\mathcal{H}_i(t) = \left\{ \frac{\partial \varphi_H(-1,\widehat{\Gamma}_H(t))}{\partial \Gamma_H(t)} \right\} \mathcal{V}^{-1}(t) \mathcal{U}(t) \mathcal{T}_i(t;\mathcal{Z})$.
Finally, by the Central Limit Theorem and similar derivations in Wang and Wang (2015), we conclude that $\sqrt{n} \left\{ \widehat{H}_{SIMEX}(t) - H_0(t) \right\}$ converges to the Gaussian process with mean zero and covariance function $E\left\{ \mathcal{H}_i(t) \mathcal{H}_i(s) \right\}$. 

$\hfill \square$

\end{document}